\documentclass[12pt]{article} 
\usepackage{amsmath}
\usepackage{amssymb} 
\usepackage{amscd}
\usepackage{amsthm}
\usepackage[utf8]{inputenc}
\usepackage{hyperref}

\newtheorem{theorem}{Theorem}

\newtheorem{proposition}[theorem]{Proposition}

\def\dim{{\mbox{dim}}}

\def\calp{{\mathcal P}} 
\def\cala{{\mathcal A}} 
\def\calb{{\mathcal B}} 
\def\cald{{\mathcal D}}

\def\calm{{\mathcal M}} 
 
\def\calo{{\mathcal O}} 
\def\calh{{\mathcal H}}

\def\calc{{\mathcal C}}

\def\calr{{\mathcal R}}

\def\fraca{{\mathfrak A}}
\def\fracb{{\mathfrak B}} 
 
 \def\fracg{{\mathfrak g}}

\def\fracr{{\mathfrak R}}

\def\bbbone{\mbox{\rm 1\hspace {-.6em} l}}

\begin{document}

\enlargethispage{3cm}

\thispagestyle{empty}
\begin{center}
{\bf LECTURES ON THE CLASSICAL MOMENT PROBLEM AND ITS NONCOMMUTATIVE GENERALIZATION}
\end{center}
   
\vspace{0.3cm}

\begin{center}
Michel DUBOIS-VIOLETTE
\footnote{Laboratoire de Physique Th\'eorique, UMR 8627, 
CNRS et Universit\'e Paris-Sud 11,
B\^atiment 210, F-91 405 Orsay Cedex\\
Michel.Dubois-Violette$@$u-psud.fr} \end{center}
 \vspace{0,5cm}


\begin{abstract}
These notes contain a presentation of the noncommutative generalization of the classical moment problem introduced in \cite{mdv:1975b} and \cite{mdv:1977}. They also contain a short summary of the classical moment problem in infinite dimension.
\end{abstract}

\vfill
 
 \newpage
\tableofcontents

\section{Introduction}\label{intro}

The aim of these lectures is to give a presentation of the noncommutative generalization of the classical moment problem introduced and studied in \cite{mdv:1975b} and \cite{mdv:1977} and to compare it with its commutative counterpart namely the classical moment problem. In order to do this we give an appropriate description of the classical moment problem and since here we do not intend to discuss notions of dimension in the noncommutative setting, our description should apply to the case of infinite dimensional spaces. This is why we give first a short summary of the relevant part of measure theory and of the classical moment problem in this context. We then describe the noncommutative generalization of the classical moment problem called the $m$-problem. In this generalization the algebra of complex polynomials is replaced by an arbitrary unital $\ast$-algebra $\fraca$ which is separated by its $C^\ast$-semi-norms, the sequence of the moment problem is replaced by a linear form on $\fraca$, the measures are replaced by positive linear forms on a $C^\ast$-algebra canonically associated with $\fraca$ and substitutes for the integration formulas of the classical moment problem are given. The connection between determination of the classical moment problem and the self-adjointness properties in the corresponding  (unbounded) representations of the polynomials algebras are generalized. A remarkable property of tensor algebras which generalizes the solubility of the one-dimensional Hamburger's moment problem is pointed out. In the case where $\fraca$ is a locally convex $\ast$-algebra we introduce a continuity condition on the solutions of the problem which generalizes the continuity condition on cylindrical measures \cite{gue-vil:1967} connected with the notion of ``scalar cylindrical concentration" \cite{sch:1973}.

Noncommutative measure theory has a very rich structure with no classical counterpart as shown by Alain Connes (occurrence of canonically associated dynamical systems) \cite{ac:1973}. We do not discuss this subject in these lectures in spite of the fact that one meets this structure in the applications of the noncommutative moment problem (the $m$-problem) to quantum fields where factors of type III$_1$ enter and where the corresponding dynamical systems should get a physical interpretation.

For the proofs of the statement concerning the noncommutative moment problem we refer to \cite{mdv:1975b} and \cite{mdv:1977}.
As explained in details in \cite{mdv:1975b} many statements there are easy consequences of powerful results of H.J. Borchers in 
\cite{bor:1962}, \cite{bor:1972} and \cite{bor:1973}.

\section{Preliminaries on $\ast$-algebras}

\subsection{Definitions}

In the following, a $\ast$-algebra $\fraca$ is an associative complex algebra $\fraca$ endowed with an antilinear involution $x\mapsto x^\ast$ such that 
\[
(xy)^\ast=y^\ast x^\ast
\]
for any $x,y \in \fraca$. An element $x$ of $\fraca$ is said to be hermitian if $x^\ast=x$. We denote by $\fraca^h$ the real subspace of all hermitian elements of $\fraca$.\\

A $C^\ast$-algebra is a $\ast$-algebra $\fracb$ which is a Banach space (i.e. a complete normed space) for a norm $x\mapsto \parallel x\parallel$ satisfying
$\parallel xy\parallel \leq \parallel x\parallel \parallel y\parallel$ and $
\parallel x^\ast x\parallel=\parallel x\parallel^2$
for any $x,y\in \fracb$. This implies $\parallel x^\ast\parallel=\parallel x\parallel$ \cite{dix:1964}, \cite{sak:1971}.\\

A $W^\ast$-algebra is a $C^\ast$-algebra $\fracr$ which is the dual Banach space ($\fracr_\ast)'$ of a Banach space $\fracr_\ast$. It can be shown that then the Banach space $\fracr_\ast$ is unique, it is called the predual of the $W^\ast$-algebra $\fracr$ \cite{sak:1971}.\\

A linear form $\phi$ on a $\ast$-algebra $\fraca$ is said to be positive if one has 
\[
\phi(x^\ast x) \geq 0
\]
for any $x\in \fraca$.

\subsection{The GNS construction}

Let $\fraca$ be a unital $\ast$-algebra. With any positive linear form $\phi$ on $\fraca$ (i.e. $\phi(x^\ast x)\geq 0$, $\forall x \in \fraca$) is associated a Hausdorff pre-Hilbert space $D_\phi$,  an element $\Omega_\phi$ of $D_\phi$ and a homomorphism $\pi_\phi$ of associative algebras with units of $\fraca$ into the algebra of endomorphisms of $D_\phi$ satisfying $D_\phi=\pi_\phi(\fraca) \Omega_\phi$, $\phi(x)=(\Omega_\phi\vert \pi_\phi(x)\Omega_\phi)$, $(\Phi \vert \pi_\phi(x)\Psi)=(\pi_\phi(x^\ast)\Phi\vert \Psi)$ for any $x\in \fraca$ and $\Phi,\Psi\in D_\phi$.  Let $\calh_\phi$ be the Hilbert space obtained by completion of $D_\phi$ ; the quadruplet $(\pi_\phi, D_\phi,\Omega_\phi, \calh_\phi)$ is unique under the above conditions up to a unitary equivalence ; it is a (generally unbounded) $\ast$-representation of $\fraca$ called {\sl the representation associated with $\phi$} \cite{pow:1971}.

\subsection{Tensor $\ast$-algebras}

Let $E$ be a real vector space and let $T_\mathbb C (E)$ be the tensor algebra over the complexified space of $E$ equipped with its structure of complex algebra with unit and the unique antilinear involution, $x\mapsto x^\ast$, for which $E(\subset T_\mathbb C(E))$ is pointwise invariant and $(xy)^\ast=y^\ast x^\ast$ $\forall x,y\in T_\mathbb C(E)$. Then $T_\mathbb C(E)$ is a $\ast$-algebra with unit which we call the tensor $\ast$-algebra over $E$. This unital $\ast$-algebra together with the canonical embedding $E\subset T_\mathbb C(E)^h$ of $E$ is characterized by the following universal property : Let $\fraca$ be a unital $\ast$-algebra, then any $\mathbb R$-linear mapping, $\alpha:E\rightarrow \fraca^h$, of $E$ into the real vector space $\fraca^h$ lifts uniquely as an homomorphism $T_\mathbb C(\alpha) : T_\mathbb C(E)\rightarrow \fraca$ of unital $\ast$-algebras.

\noindent \underbar{Remark}. Since $x\mapsto x^\ast$ is canonically an anti-isomorphism of $T_\mathbb C(E)$ onto the opposite $\ast$-algebra, the above property is equivalent to the following one : Any $\alpha:E\rightarrow \fraca^h$ as above lifts uniquely as an anti-homomorphism $\tilde T_\mathbb C(\alpha):T_\mathbb C(E)\rightarrow \fraca$. This remark is relevant for quantum field theory \cite{str-wig:1964} because when $E$ is the space of real test functions, the Borchers field algebra \cite{bor:1962} is just the completion of $T_\mathbb C(E)$ for a suitable topology and it is known that some space-time symmetries have to be represented there by automorphisms and some others by anti-automorphisms (e.g. TCP) of the Borchers algebra.

\section{Cylindrical measures and the classical moment problem}

\subsection{Polynomials and cylindrical functions} . 

Let $E$ be a real vector space with algebraic dual $E^\ast$. Suppose that instead of working with unital $\ast$-algebras we are only interested in commutative unital $\ast$-algebras. Then the analog of $T_\mathbb C(E)$ is the symmetric $\ast$-algebra over $E$ denoted bt $S_\mathbb C(E)$. This is the complex symmetric algebra over the complexified vector space of $E$ equipped with the unique anti-linear involution leaving $E$ pointwise invariant and such that it is a commutative $\ast$-algebra with unit. $S_\mathbb C(E)$ is also characterized by a universal  property. Any $\mathbb R$-linear mapping $\alpha:E\rightarrow \fraca^h$ of $E$ into the real vector space $\fraca^h$ of the hermitian elements of a commutative $\ast$-algebra with unit $\fraca$ lifts uniquely as an homomorphism $S_\mathbb C(\alpha):S_\mathbb C(E)\rightarrow \fraca$ of commutative $\ast$-algebras with units. Let $S_\mathbb C(E)^\wedge$ denote the set of (characters of $S_\mathbb C(E)$) all the $\ast$-homomorphisms $\chi$ of $S_\mathbb C(E)$ into $\mathbb C$ mapping the unit of $S_\mathbb C(E)$ onto $1\in \mathbb C$ $(\chi(\bbbone)=1)$. The restriction to $E\subset S_\mathbb C(E)$ maps $S_\mathbb C(E)^\wedge$ into $E^\ast$ and it follows from the above universal properly applied to the case $\fraca=\mathbb C$ (so $\fraca^h=\mathbb R$) that it is a bijection $(\chi=S_\mathbb C(\chi) \restriction E)$ of $S_\mathbb C(E)^\wedge$ onto $E^\ast$. Let $p$ be an element of $S_\mathbb C(E)$ and $\xi$ an element of $E^\ast$ ; the value at $p$ of $S_\mathbb C(\xi)\in S_\mathbb C(E)^\wedge$ will simply be denoted by $p(\xi)$. Let $E^\ast_\sigma$ be $E^\ast$ equipped with the weak topology $\sigma(E^\ast,E)$. Then $\xi\mapsto p(\xi)$ is for each $p\in S_\mathbb C(E)$, a continuous function on $E^\ast_\sigma$; we call these functions {\sl polynomial functions on $E^\ast_\sigma$}. These functions form a $\ast$-subalgebra  with unit of the algebra $\mathbb C^{E^\ast}$ of all complex functions on $E^\ast$ which is isomorphic to $S_\mathbb C(E)$ (under $p\mapsto (\xi\mapsto p(\xi))$). For any $p\in S_\mathbb C(E)$ there is a finite family $h_1,\cdots, h_n$ in $E$ and a polynomial function $P$ on $\mathbb R^n$ for which $p(\xi)=P(\langle h_1,\xi\rangle,\cdots, \langle h_n,\xi))$ ($\forall \xi \in E^\ast$). More generally a cylindrical function on $E^\ast_\sigma$ is a function on $E^\ast$ of the form $\xi\mapsto f(\langle h_1,\xi\rangle, \cdots, \langle h_n,\xi\rangle)$ for some finite family $h_1,\cdots, h_n$ in $E$ and some complex function $f$ on $\mathbb R^n$. These functions also form a $\ast$-subalgebra with unit of $\mathbb C^{E^\ast}$. Let $h_1,\cdots,h_n$ be a finite family in $E$, we denote by $\calc_{(0)}(h_1,\cdots, h_n)$ the set of cylindrical functions $\xi\mapsto f(\langle h_1,\xi\rangle, \cdots, \langle h_n,\xi\rangle)$ when $f$ runs over the $C^\ast$-algebra $\calc_{(0)}(\mathbb R^n)$ of complex continuous functions vanishing at infinity on $\mathbb R^n$. This is a $C^\ast$-subalgebra of the $C^\ast$-algebra $\calc^b(E^\ast_\sigma)$ of complex continuous bounded function on $E^\ast_\sigma$.  Let $\calb(S_\mathbb C(E),E)$ be the $C^\ast$-subalgebra of $\calc^b(E^\ast_\sigma)$ generated by $\bigcup_{h\in E}\calc_{(0)}(h)$; it contains $\cup\calc_{(0)}(h_1,\cdots,h_n)$, where the union is taken over the finite families in $E$, as a dense $\ast$-subalgebra. Let us set $f(h_1,\cdots, h_n)(\xi)=f(\langle h_1,\xi\rangle),\cdots,\langle h_n,\xi\rangle)$.

\subsection{Cylindrical measures}

Here $E$ is again a real vector space and we use the above notations. We say that a positive linear form $\omega$ on $\calb(S_\mathbb C(E),E)$ has {\sl the property} (C) if it satisfies the following condition :
\begin{equation}\tag{C}
\forall h\in E, \>\>\> \parallel \omega\restriction \calc_{(0)}(h)\parallel =\parallel \omega\parallel
\end{equation}
$(\omega\restriction \calc_{(0)}(h)$ denotes the restriction of $\omega$ to $\calc_{(0)}(h)\subset \calb(S_\mathbb C(E),E))$. It can be shown that the property (C) for $\omega$ is equivalent to the following (a priori stronger) property (C$^\prime$) : 
\begin{equation}\tag{C$^\prime$}
\parallel\omega \restriction \calc_{(0)}(h_1,\cdots, h_n)\parallel = \parallel \omega\parallel
\end{equation}
for any finite family $(h_1,\cdots, h_n)$ in $E$. Thus $f\mapsto f(h_1,\cdots, h_n)$ is a $\ast$-homomorphism of $\calc_{(0)}(\mathbb R^n)$ in $\calb(S_\mathbb C(E),E)$ and therefore it follows that $f\mapsto \omega(f(h_1,\cdots, h_n))$ is a positive linear form on $\calc_{(0)} (\mathbb R^n)$ for any positive linear form $\omega$ on $\calb(S_{\mathbb C}(E),E)$. By the Riesz theorem we have 
\[
\omega (f(h_1,\cdots, h_n))=\int f d\mu_{h_1,\cdots,h_n}
\]
for a unique positive bounded measure $\mu_{h_1,\cdots, h_n}$ on $\mathbb R^n$. If furthermore, $\omega$ has the property (C), this system of measures is coherent in the following sense : If  $f$ and $f'$ are bounded Borel functions on $\mathbb R^n$ and $\mathbb R^{n'}$ such that we have $f(\langle h_1,\xi\rangle,\cdots, \langle h_n,\xi\rangle)=f'(\langle h'_1,\xi\rangle, \cdots, \langle h'_{n'},\xi\rangle)$ for any $\xi \in E^\ast$ then we also have $\int f d\mu_{h_1,\cdots, h_n}=\int f'd\mu_{h'_1,\cdots, h'_{n'}}$, (in particular $\int d\mu_{h_1,\cdots, h_n}=\int d\mu_{h'_1,\cdots, h'_n}=\parallel \omega \parallel$). Such a coherent family of positive bounded measure labelled by the finite families in $E$ is called a cylindrical measure on $E^\ast_\sigma$. Conversely, given a cylindrical measure on $E^\ast_\sigma$, we define a positive linear form $\omega$ on $\cup\calc_{(0)}(h_1,\cdots, h_n)$ by $\omega(f(h_1,\cdots, h_n))=\int f d\mu_{h_1,\cdots,h_n}$ which has norm $\parallel \omega\parallel=\int d\mu_{h_1,\cdots,h_n}<\infty$ and, therefore, extends uniquely into a positive linear form on $\calb(S_{\mathbb C}(E),E)$ which obviously has property (C). If $\mu$ is a positive bounded regular measure on $E^\ast_\sigma$ (i.e. a Radon measure in the sense of \cite{sch:1973}), then $\omega(f)=\int_{E^n} fd\mu$ defines (for $f\in \calb(S_\mathbb C(E),E)$) a positive linear form $\omega$ on $\calb(S_\mathbb C(E),E)$ which has the property (C). Therefore, with a bounded positive regular measure $\mu$ on $E^\ast_\sigma$ (and a fortiori on $E'_\sigma$ for any dense subspace $E'$ of $E^\ast_\sigma$) is associated (in an injective way) a cylindrical measure on $E^\ast_\sigma$. If $E$ is finite dimensional the converse is also true ; but it is wrong if $E$ has an uncountable basis and the set of bounded positive regular measures on $E^\ast_\sigma$ is (canonically), in this case, a strict subset of the set of cylindrical measures on $E^\ast_\sigma$.

\subsection{Continuity condition and Minlos theorem}

Let $E_\theta$ be a real vector space $E$ equipped with a Hausdorff locally convex topology $\theta$. The weak dual $E'_\sigma$ of $E_\theta$ is canonically a dense topological vector subspace of $E^\ast_\sigma$. It follows that the restriction $f\mapsto f\restriction E'$ defines an isomorphism between $\calb(S_\mathbb C(E),E)$ and the corresponding $C^\ast$-algebra of continuous functions on $E'_\sigma$ and that, more generally, there are no distinctions between continuous functions on $E'_\sigma$ and continuous functions on $E^\ast_\sigma$. In a dual manner, there are no distinctions between cylindrical measures on $E'_\sigma$ and on $E^\ast_\sigma$ (although measures on $E'_\sigma$ and $E^\ast_\sigma$ are distinct). There is, however, a very natural way to select by continuity a subset of cylindrical measures associated with the topology $\theta$. Namely we say that a cylindrical measure (on $E^\ast_\sigma$) satisfies the $\theta$-{\sl continuity condition} if, for any integer $n$ and family $f_1,\cdots, f_n$ in $\calc_{(0)}(\mathbb R)$, the corresponding positive linear form $\omega$ on $\calb(S_\mathbb C(E),E)$ is such that $(h_1,\cdots, h_n)\mapsto \omega(f_1(h_1)\cdots f_n(h_n))$ is a continuous function on $E^n_\theta$. The theorem of Minlos \cite{min:1963}, \cite{gue-vil:1967}, 
\cite{sch:1973}, \cite{cho:1969} consists in the following : If  $E_\theta$ is a nuclear space then a cylindrical measure on $E^\ast_\sigma$ which satisfies the $\theta$-continuity condition is a positive bounded regular measure on $E'_\sigma$, (notice that if $E_\theta$ is a barreled space, every bounded positive regular measure on $E'_\sigma$ satisfies this condition).

\subsection{Integration of cylindrical functions}   

Let $E$ be a real vector space and $(\mu_{h_1,\cdots, h_n})$ be a cylindrical measure on $E^\ast_\sigma$, we denote by $\omega$ the corresponding positive linear form on $\calb(S_\mathbb C(E),E)$. Let $f$ be a Borel function on $\mathbb R^n$ which is in $L^1(d\mu_{h_1,\cdots, h_n})$. Then the integral $\int f d\mu_{h_1,\cdots,h_n}$ does only depend on the function $f(h_1,\cdots,h_n)$ on $E^\ast$ ;  we denote it by $\bar\omega(f(h_1,\cdots, h_n))$ and call it {\sl the integral of the cylindrical function} $f(h_1,\cdots,h_n)$. If every polynomial $p\in S_\mathbb C(E)$ on $E^\ast_\sigma$ is integrable (i.e. all the $\mu_{h_1,\cdots,h_n}$ are rapidly decreasing measures) we say that $(\mu_{h_1,\cdots,h_n})$ is rapidly decreasing. In this case $\bar \omega$ defines a positive linear form on the $\ast$-algebra $\calp(E^\ast_\sigma)$ of continuous polynomially bounded cylindrical functions on $E^\ast_\sigma$. Conversely, if $\psi$ is  a positive linear form on $\calp(E^\ast_\sigma)$, its restriction to $\calb(S_\mathbb C(E),E)$ has the property (C) and the corresponding cylindrical measure is rapidly decreasing and $\psi(f)=\bar\omega(f)$ $(\forall f\in \calp(E^\ast_\sigma))$.

\subsection{The classical moment problem}

Let us use the above notations and suppose that $(\mu_{h_1,\cdots, h_n})$ is a rapidly decreasing cylindrical measure on $E^\ast_\sigma$ ; then ``{\sl its moments}" 
\[
\int t_1\cdots t_n d\mu_{h_1,\cdots,h_n}(t_1,\cdots, t_n)=\bar \omega(h_1,\cdots h_n)
\]
exist $(\forall h_1\cdots, h_n$ in $E$). These moments correspond to a unique linear form $\phi$ on $S_\mathbb C(E)$ $(\phi(p)=\bar \omega(p)$, $\forall p\in S_\mathbb C(E))$ which is positive on the positive valued polynomials on $E^\ast_\sigma$ ; we say that $\phi$ is a {\sl strongly positive} linear form on $S_\mathbb C(E)$. Such a strongly positive linear form is a positive linear form on the $\ast$-algebra $S_\mathbb C(E)$ but if $\dim(E)\geq 2$, there are positive linear forms on $S_\mathbb C(E)$ which are not strongly positive linear forms. In the case $\dim(E)=1$, both concepts coincide but in this case $S_\mathbb C(E)=T_\mathbb C(E)$ and this coincidence will appear as a specific case of a general result on the algebras $T_\mathbb C(E)$ (see below). In the {\sl classical moment problem}, one starts from a linear form $\phi$ on $S_\mathbb C(E)$ and one asks whether the $\phi(h_1\cdots h_n)$ are the moments of some cylindrical measure (rapidly decreasing) on $E^\ast_\sigma$ (notice that if this is the case for a measure $\mu$ on $E^\ast_\sigma$ this means  $\phi(h_1\cdots h_n)=\int \langle h_1,\xi\rangle \cdots \langle h_n,\xi\rangle d\mu(\xi)$. In order that this should happen, $\phi$ must be strongly positive, and it turns out that it is sufficient. To see this, we notice that, by a direct application of the Hahn-Banach theorem, a strongly positive linear form on $S_\mathbb C(E)$ has positive extensions to the cylindrical continuous polynomially bounded functions on $E^\ast_\sigma$ and that (see 3.4) these extensions are canonically rapidly decreasing cylindrical measures on $E^\ast_\sigma$. We call these cylindrical measures {\sl solution of the moment problem for} $\phi$. If $\theta$ is a Hausdorff locally convex topology on $E$ and if for each integer $n, (h_1,\cdots, h_n)\mapsto \phi(h_1\cdots h_n)$ is continuous on $E^n_\theta$ then any solution of the moment problem satisfies the $\theta$-continuity condition. It follows (by Minlos theorem) that if $\theta$ is a nuclear topology and if $\phi$ has the above continuity property then any solution of the moment problem for $\phi$  is a bounded positive regular measure on the weak dual $E'_\sigma$ of $E_\theta$ \cite{heg:1975}, \cite{sho-tam:1963}. If $\phi$ is such that there is a unique solution of the moment problem, we say that the moment problem for $\phi$ is {\sl determined} ; there is a classical connection between determination and self-adjointness properties in the GNS representation associated with $\phi$ \cite{sho-tam:1963}, \cite{akh:1965} which will be generalized in the noncommutative case (see in Section 5).

\section{Noncommutative generalization of the classical moment problem : The $m$-problem}

In this section we shall generalize the concepts associated with the pair $(S_\mathbb C(E),E)$ described in the introduction to pairs $(\fraca,E)$ where $\fraca$ is a (noncommutative) $\ast$-algebra with unit and $E$ is a real vector subspace of $\fraca^h$ which is {\sl generating} for $\fraca$ (as a unital $\ast$-algebra, i.e. an $\ast$-subalgebra of $\fraca$ which contains the unit and $E$ must be identical with $\fraca$). In particular we shall give a (noncommutative) generalization of cylindrical measures and of the classical moment problem (the $m$-problem). Special attention will be devoted to the pair $(T_\mathbb C(E),E)$ because it follows from the universal property of $T_\mathbb C(E)$ that there is a unique surjective homomorphism $T_\mathbb C(E)\rightarrow \fraca$ of unital  $\ast$-algebras which induces the identity of $E$ onto itself (notice that this is the very reason why these tensor algebras enter into the formulation of quantum field theory \cite{bor:1962}); so $\fraca$ is a quotient of $T_\mathbb C(E)$.

\subsection{$C^\ast$-semi-norms}

Let $\fraca$ be a $\ast$-algebra with unit. A $C^\ast$-{\sl semi-norm} on $\fraca$ will be a semi-norm $q$ on $\fraca$ which satisfies :
$q(xy)\leq q(x)q(y), q(x^\ast)=q(x),q(x^\ast x)=q(x)^2(\forall x,y\in \fraca)$ and $q(\bbbone)=1$. In other words $q$ is such that the completion $\hat\fraca_q$ of $\fraca/q^{-1}(0)$ for induced norm is canonically a $C^\ast$-algebra with unit \cite{dix:1964}. Notice that any $C^\ast$-semi-norm $q$ on $S_\mathbb C(E)$ is of the form $p\mapsto q(p)=\sup\{\vert p(\xi)\vert \vert \xi\in K_q\}$ where 
$K_q=\{\xi\in E^\ast \parallel p(\xi)\vert \leq q(p), \forall p\in S_\mathbb C(E)\}$ is a compact subset of $E^\ast_\sigma$ (every closed bounded subset of $E^\ast_\sigma$ is compact \cite{sch:1971}),  and that conversely $p\mapsto \sup\{\vert p(\xi)\vert \vert \xi\in K\}$ is a $C^\ast$-semi-norm on $S_\mathbb C(E)$ for any compact subset $K$ of $E^\ast_\sigma$. It follows (using the Stone-Weierstrass theorem) that the completion of $S_\mathbb C(E)$ for the locally convex topology generated by its $C^\ast$-semi-norms is canonically the $\ast$-algebra of all functions on $E^\ast_\sigma$ which are continuous on the compact subsets of $E^\ast_\sigma$. From this algebra, it is then not hard to extract all the concepts which enter into the formulation of the classical moment problem. It is the noncommutative generalization of this algebraic construction that we shall now describe (for the proofs, we refer to \cite{mdv:1975b}). In this section we shall consider unital $\ast$-algebras which are separated by their $C^\ast$-semi-norms so the following result is worth noticing.

\subsection{A remarkable property of tensor algebras}
\begin{theorem}
Let $E$ be a real vector space and let $T_\mathbb C(E)^+$ denote the convex hull in $T_\mathbb C(E)$ of $\{x^\ast x\vert x\in T_\mathbb C(E)\}$. Then $T_\mathbb C(E)^+$ is a convex salient cone in $T_\mathbb C(E)$ which is closed for the locally convex topology generated by the $C^\ast$-semi-norms on $T_\mathbb C(E)$.
\end{theorem}
For the proof see \cite{mdv:1975b}, Theorem 2.\\

This theorem implies, in particular, that the locally convex topology generated by the $C^\ast$-semi-norms on $T_\mathbb C(E)$ is separated (since the closure of $\{0\}$ is contained in $T_\mathbb C(E)^+ \cap (-T_\mathbb C(E)^+)=\{0\})$. But it implies much more ; for instance this is wrong if $T_\mathbb C(E)$ is replaced by $S_\mathbb C(E)$ with $\dim(E)\geq 2$ because there, the closure of $S_\mathbb C(E)^+$ for the locally convex  topology generated by the $C^\ast$-semi-norms is (see above) the set of all the positive valued polynomials on $E^\ast_\sigma$ which is strictly bigger than $S_\mathbb C(E)^+$ = (finite sums of squares of absolute values).

\subsection{Properties of the completions}
\begin{proposition}
Let $\fraca$ be a $\ast$-algebra with unit such that the locally convex topology generated by its $C^\ast$-semi-norms is separated and let $\cala$ be the topological $\ast$-algebra with unit obtained by completion of $\fraca$ for this topology. We denote, as usual, by $\cala^h$ the $\mathbb R$-vector space of all hermitian elements of $\cala$  and by $\cala^+$ the convex hull in $\cala$ of $\{x^\ast x\vert x\in \cala\}$ and for $x\in \cala$ we let $Sp(x)$ be the spectrum of $x$ in $\cala$ (i.e. $Sp(x)=\{\lambda\in \mathbb C\vert (x-\lambda\bbbone)$ has no inverse in $\cala\}$). Then we have the following .\\
a) $x \in \cala^h$ is equivalent to $Sp(x)\subset \mathbb R$.\\
b) $\cala^+$ is a closed convex salient cone in $\cala$ and we have :
\[
\cala^+=\{x\in \cala\vert Sp(x)\subset \mathbb R^+\}= \{h^2\vert h\in \cala^h\}
\]
c) Every $h\in \cala^h$ has a unique decomposition $h=h^+-h^-$ with $h^+,h^-\in \cala^+$ and $h^+h^-=0$\\
d) Every $h\in \cala^h$ determines a unique homomorphism $f\mapsto f(h)$ of $\ast$-algebras with units from the $\ast$-algebra $\calp(\mathbb R)$ of all complex continuous polynomially bounded functions on $\mathbb R$ into the $\ast$-algebra $\cala$ such that $Id_\mathbb R(h)=h$ (where $Id_\mathbb R$ denote the identity mapping $t\mapsto t$, of $\mathbb R$ onto itself). This homomorphism is (automatically) continuous if $\calp(\mathbb R)$ is equipped with the topology of compact convergence on $\mathbb R$ and $f(h)=0$ is equivalent to $f\restriction Sp(h)=0$.
\end{proposition}
Notice that 
\[
\calb_\infty=\{x\in \cala\vert \parallel x\parallel=\sup\{q(x)\vert q\in C^\ast-\text{semi-norms on}\ \fraca\}\leq \infty\}
\]
 is canonically a $C^\ast$-algebra and that if $h\in \cala^h$ and if $f$ is a continuous bounded function on $\mathbb R$, then $f(h)\in \calb_\infty$.

\subsection{The $m$-problem}

Let $\fraca$ be as above and $E$ be a real vector subspace of $\fraca^h$ which generates $\fraca$ as $\ast$-algebra with unit. We define the $C^\ast$-{\sl algebra} $\calb(\fraca,E)$ associated with the pair ($\fraca,E)$ to be the $C^\ast$-subalgebra of $\calb_\infty$ generated by $\{f(h)\vert f\in \calc_{(0)}(\mathbb R), h\in E\}$. $\calb(\fraca,E)$ is a dense $\ast$-subalgebra of $\cala$ (which generally does not contain the unit) and, for $h\in E$, $\calc_{(0)}(h)=\{f(h)\vert f \in \calc_{(0)}(\mathbb R)\}$ is a $C^\ast$-subalgebra of $\calb(\fraca,E)$.
The notations are coherent with the notations introduced in 3.1 (for the case $\fraca=S_\mathbb C(E)$) and we say that a positive linear form $\omega$ on $\calb(\fraca,E)$ has {\sl the property} (C) if $\parallel \omega\restriction \calc_{(0)}(h)\parallel=\parallel \omega\parallel$ for any $h$ in $E$. This is the noncommutative generalization of the notion of cylindrical measure we need to formulate our generalization of the moment problem on $(\fraca,E)$. We say that a linear form $\phi$ on $\fraca$ is {\sl strongly positive} if it is positive valued on the closure $\fraca \cap \cala^+$ of $\fraca^+$ in $\fraca$ equipped with the locally convex topology generated by its $C^\ast$-semi-norms. Let $\phi$ be a linear form on $\fraca$; we say that a positive linear form $\omega$ on $\calb(\fraca,E)$ is a solution of the $m$-{\sl problem for} $\phi$ on $(\fraca, E)$ or of the $m(E)$-{\sl problem for} $\phi$ if $\omega$ has a positive extension $\bar\omega$ on some ordered subspace $\calm$ of $\cala$ (equipped with its positive cone $\cala^+$) which contains $\fraca$ and such that $\phi=\bar\omega\restriction \fraca$, i.e. $\phi(h_1\cdots h_n)=\bar\omega(h_1\cdots h_n)$ for $h_1,\cdots, h_n\in E$.

\subsection{Solubility condition for the $m$-problem}
\begin{theorem}
Let $(\fraca,E)$ be as above, and $\phi$ be a linear form on $\fraca$. Then the $m(E)$-problem for $\phi$ is soluble if and only if $\phi$ is strongly positive. The set $\fracg_\phi(E)$ of all solutions of the $m(E)$-problem for $\phi$ is convex and compact in the weak dual of $\calb(\fraca,E)$; if $\phi'\not=\phi$ then $\fracg_{\phi'}(E)\cap \fracg_{\phi}(E)$ is empty. If $\omega\in \fracg_\phi(E)$, then $\omega$ has the property (C) and, more generally we have : $\phi(h^n)=\int t^n d\mu_{h,\omega}(t)$ for any $h\in E$ and integer $n\geq 0$; where the measure $\mu_{h,\omega}$ on $\mathbb R$ is defined by $\omega(f(h))=\int f(t) d\mu_{h,\omega}(t)$, $\forall f\in \calc_{(0)}(\mathbb R)$.
\end{theorem}
We notice that Theorem 1 means that every positive linear form on $T_\mathbb C(E)$ is strongly positive so the last theorem implies that the $m$-problem is always soluble for a positive linear form on $T_\mathbb C(E)$. This generalizes the solubility of one-dimensional Hamburger problem for positive moment sequences since $\mathbb C[X]\simeq S_\mathbb D(\mathbb R)=T_\mathbb C(\mathbb R)$ and since it is clear from the introduction that when $\fraca=S_\mathbb C(E)$ then the $m(E)$-problem reduces to the classical moment problem.\\
If $\fracg_\phi(E)$ has exactly one element, we say that the $m(E)$-problem for $\phi$ is {\sl determined}.

\subsection{Convergence}
\begin{proposition}
Let $\phi_\alpha$ be a net of strongly positive linear forms on $\fraca$ which converges weakly to $\phi$ and let $\omega_\alpha\in\fracg_{\phi_\alpha}(E)$ for each $\alpha$. Then, $\phi$ is strongly positive and any weak limit of a subnet of $\omega_\alpha$ is in $\fracg_\phi(E)$.
\end{proposition}
If $\phi_\alpha(\bbbone)$ is bounded, then $\parallel \omega_\alpha\parallel =\phi_\alpha(\bbbone)$ is bounded so the $\omega_\alpha$ belong to a weakly compact set of positive linear forms on $\calb(\fraca, E)$ and, therefore in this case there are weakly convergent subnets of $\omega_\alpha$.

\subsection{Homomorphisms}

If $(\fraca_i,E_i)\>\> i=1,2$ are pairs of $\ast$-algebras with units $\fraca_i$ and generating subspace $E_i\subset \fraca_i^h$, we define a {\sl morphism} $\alpha : (\fraca_1,E_1)\rightarrow (\fraca_2,E_2)$ to be a homomorphism $\alpha$ of $\ast$-algebras with unit from $\fraca_1$ into $\fraca_2$ such that $\alpha(E_1)\subset E_2$. It is then not hard to see that for any morphism $\alpha:(\fraca_1,E_1)\rightarrow (\fraca_2,E_2)$ there is a unique $\ast$-homomorphism $\calb(\alpha):\calb(\fraca_1,E_1)\rightarrow \calb(\fraca_2,E_2)$ for which we have $\calb(\alpha)(f(h))=f(\alpha(h))$, $\forall h\in E_\alpha$ and $f\in \calc_{(0)}(\mathbb R)$. Furthermore, this correspondence is functorial.

 If $E$ is a real vector space, we denote the $C^\ast$-algebra $\calb(T_\mathbb C(E),E)$ by $\calb_0(E)$. If $A:E_1\rightarrow E_2$ is a linear mapping from a real vector space $E_1$ into another one $E_2$, then using the universal property of $T_\mathbb C(E)$ we obtain a unique $\ast$-homomorphism $\calb_0(A):\calb_0(E_1)\rightarrow \calb_0(E_2)$ which satisfied $\calb_0(A)(f(h))=f(Ah)$, $\forall h\in E_1$ and $f\in \calc_{(0)}(\mathbb R)$ (and $\calb_0$ is a covariant fonctor from the category of real vector spaces into the category of $C^\ast$-algebras). However, in view of the remark made at the end of 2.3, there is also, for $A$ as above, a unique anti $\ast$-homomorphism $\tilde\calb_0(A):\calb_0(E_1)\rightarrow \calb_0(E_2)$ satisfying $\tilde\calb(A)(f(h))=f^\ast(Ah)$, $\forall h\in E_1$ and $f\in \calc_{(0)}(\mathbb R)$.

\section{Representations, self-ajdointness and determination}

In this section, $(\fraca, E)$ is again a pair which consists of a $\ast$-algebra with unit $\fraca$ and a real vector subspace $E$ of $\fraca^h$ which generates $\fraca$ as $\ast$-algebra with unit and it is assumed that $\fraca$ is separated by its $C^\ast$-semi-norms.

\subsection{Self-adjointness and determination}
\begin{theorem}
Let $\phi$ be a strongly positive linear form on $\fraca$ and let $\omega$ be a solution of the $m(E)$-problem for $\phi$ ($\omega\in \fracg_\phi(E)$). let $(\pi_\phi,D_\phi, \Omega_\phi,\calh_\phi)$ be the representation of $\fraca$ associated with $\phi$ and $(\pi_\omega,\Omega_\omega,\calh_\omega)$ be the representation of $\calb(\fraca,E)$ associated with $\omega$. Then, $\calh_\phi$ is canonically a Hilbert subspace of $\calh_\omega$ with $\Omega_\phi=\Omega_\omega$. If, for every $h\in E$, $\pi_\phi(h)$ is essentially self-adjoint (on $D_\phi$) then $\calh_\phi=\calh_\omega$ and $\pi_\omega(f(h))=f(\overline{\pi_\phi(h)})$ for $f\in \calc_0(\mathbb R)$ and  $h\in E$; so in this case, the $m(E)$-problem for $\phi$ is determined (i.e. $\omega$ is unique).
\end{theorem}
This generalizes the connection between self-adjointness and determination in the classical moment problem. For the proof, see reference \cite{mdv:1977}.

\subsection{Normality and topology}.

Let $(\fraca,E)$ be as above and let $\theta$ be a locally convex topology on $E$. We say that a positive linear form $\omega$ on $\calb(\fraca, E)$ satisfies the $\theta$-{\sl continuity condition} if, for every integer $n$ and $f_1,\cdots,f_n\in \calc_{(0)}(\mathbb R)$, $(h_1,\cdots, h_n)\mapsto \omega(f_1(h_1)\cdots f_n(h_n))$ is a continuous function on $E_\theta^n$ ($E_\theta$ denotes $E$ equipped with $\theta$). We denote by $\calr^+_\ast(\fraca,E_\theta)$ the set of all the positive linear forms on $\calb(\fraca, E_\theta)$ having the property (C) (see in 4.4) and satisfying the $\theta$-continuity condition. Using standard arguments on uniform convergence, it is easily seen that $\calr^+_\ast(\fraca,E_\theta)$ is a norm closed invariant (i.e. $\omega(\bullet)\in \calr^+_\ast (\fraca, E_\theta)\Rightarrow \omega(x^\ast(\bullet)x)\in \calr^+_\ast(\fraca,E_\theta)$ $\forall x\in \calb(\fraca,E))$ convex cone of positive linear forms on $\calb(\fraca,E_\theta)$ in other words it is a {\sl folium} \cite{haa-kad-kas:1970}. So $\calr^+_\ast(\fraca, E_\theta)$ is the cone of positive normal linear form of a $W^\ast$-algebra $\calr(\fraca, E_\theta)$ \cite{sak:1971}.

\subsection{Some related results}
\begin{theorem}
Let $(\fraca, E,\theta)$ be as above and let $\phi$ be a strongly positive linear form on $\fraca$ such that, for each integer $n$, $(h_1,\cdots,h_n)\mapsto \phi(h_1\cdots h_n)$ is a continuous function on $E^n_\theta$ and such that $\pi_\phi(h)$ is essentially self-adjoint for any $h$ in $E$. Then the unique solution of the $m(E)$-problem for $\phi$ is in $\calr^+_\ast(\fraca, E_\theta)$.
\end{theorem}
This theorem is a consequence of the following classical results on strong resolvent convergence \cite{kat:1976}, 
\cite{ree-sim:1972}. If $A_\alpha$ is a net of self-adjoint operators (with domains dom ($A_\alpha))$, if $A$ is another self-adjoint operator and if there is a core $D$ for $A$ in $\cap_\alpha$ dom $(A_\alpha)$ such lim $(A_\alpha\Phi)=A\Phi$, $\forall\Phi\in D$ then $A_\alpha$ converges to $A$ in the sense of strong resolvent convergence which turns out to be equivalent with strong convergence of $f(A_\alpha)$ to $f(A)$ for each $f\in \calc_{(0)}(\mathbb R)$ (and even for each continuous bounded function $f$ on $\mathbb R$).

The last result shows that the continuity condition is a good one to lock the $m(E)$-problem with topologies on $E$.

There is a canonical $\ast$-homomorphism of $\calb(\fraca,E)$ onto a weakly \cite{sak:1971} dense $C^\ast$-subalgebra of $\calr(\fraca,E_\theta)$ which is injective whenever the $C^\ast$-semi-norms on $\fraca$ having $\theta$-continuous restrictions to $E$ separate $\fraca$. This is the case if $\fraca=T_\mathbb C(E)$ for any Hausdorff locally convex topology $\theta$ on $E$ and the following theorem of Borchers \cite{bor:1973}, \cite{mdv:1975a} shows that, in this case, we have much more.

\begin{theorem}
Let $E$ be a real locally convex vector space with complexified $E_\mathbb C=E\oplus iE$ and let $T^{(\varepsilon)}_\mathbb C(E)$ be the locally convex direct sum of the n$^{th}$ $\varepsilon$-tensor powers 
\cite{tre:1967} $\otimes^n_\varepsilon E_\mathbb C$ of $E_\mathbb C$; $T^{(\varepsilon)}_\mathbb C(E)=\oplus_{n\geq 0} (\otimes^n_\varepsilon E_\mathbb C), (\otimes^0 E_\mathbb C=\mathbb C)$. Then for each integer $N$, the continuous $C^\ast$-semi-norms on $T^{(\varepsilon)}_\mathbb C(E)$ induce on $T^N_\mathbb C(E)=\oplus^{n=N}_{n=0}(\otimes^n E_\mathbb C)$ a locally convex topology which coincides with its topology as subspace of $T^{(\varepsilon)}_\mathbb C(E)$.
\end{theorem}

It is worth noticing here that as shown in \cite{alc-yng:1988} and \cite{yng:1988} a similar result holds for certain quadratic algebras which are ``partially commutative" quotients of tensor algebras. This includes in particular the quotient of the tensor algebra over the space of test functions by the ``locality ideal" \cite{bor:1962}.

If $E$ is a Hausdorff locally convex real vector space $\calr(T_\mathbb C(E),E)$ (resp. $\calr^+_\ast(T_\mathbb C(E),E))$ will simply be denoted by $\calr_0(E)$ (resp. $\calr^+_{0\ast}(E))$. $\calr_0(E)$ contains $\calb_0(E)$ as a weakly dense $C^\ast$-subalgebra (i.e. dense for the weak topology $\sigma(\calr_0(E)$, $\calr_{0\ast}(E))$ where $\calr_{0\ast}(E)$ is the predual of $\calr_0(E))$ and, if $A:E_1\rightarrow E_2$ is a continuous linear mapping of the locally convex real vector space $E_1$ into another one $E_2$, then $\calb_0(A)$ (resp. $\tilde \calb_0(A))$ extends itself canonically as a $W^\ast$-homomorphism $\calr_0(A):\calr_0(E_1)\rightarrow \calr_0(E_2)$ (resp. $W^\ast$-anti-homomorphism $\tilde\calr_0(A):\calr_0(E_1)\rightarrow \calr_0(E_2))$. Thus $\calr_0$ is a functor from the category of Hausdorff locally convex real vector spaces and continuous linear mappings in the category of $W^\ast$-algebras and $W^\ast$-homomorphisms, (notice that $\calr_0(A_2\circ A_1)=\calr_0(A_2)\circ \calr_0(A_1)$).

\begin{theorem}
Let $E$ be a Hausdorff locally convex real vector space and $F$ be a subspace of $E$ (equipped with the induced topology). Then $\calr_0(F)$ is canonically a $W^\ast$-subalgebra of $\calr_0(E)$ which contains the unit of $\calr_0(E)$. If $G$ is another subspace of $E$ then $\calr_0(F+G)$ is generated, as $W^\ast$-subalgebra of $\calr_0(E)$,  by $\calr_0(F)\cup\calr_0(G)$. If $\hat E$ is the completion of $E$, we have : $\calr_0(\hat E)=\calr_0(E)$. The above canonical identifications are natural with respect to the functional calculus ; i.e. if $h\in F\subset E$ and $f\in\calc_{(0)}(\mathbb R)$. then $f(h)\in \calr_0(F)$ is identified with $f(h)\in \calr_0(E)$.
\end{theorem}

\subsection{Examples}
Let $Q$ be the $\ast$-algebra generated by the Schr\"odinger representation of Heisenberg canonical commutation relations $[q,p]= i\bbbone$ and $R(Q)$ be the von Neumann algebra generated by the self-adjoint operators $p$ and $q$. Let $\Omega$ be the ground state of the harmonic oscillator ; then $(t_1,t_2)\mapsto t_1p+t_2q$ is linear from $\mathbb R^2$ into $Q^h$, so there is a unique homomorphism of unital $\ast$-algebra of $T_\mathbb C(\mathbb R^2)$ onto $Q$, $\pi$ for which $\pi(t_1,t_2)=t_1p + t_2q$ and the state $x\mapsto (\Omega\vert \pi(x)\Omega)= \phi(x)$ has the property that $\pi_\phi(t_1,t_2)=\pi(t_1,t_2)\restriction D_\phi$ are essentially self-ajoint. So there is a unique solution of the $m(\mathbb R^2)$-problem for $\phi$ which (by construction) satisfies the continuity condition. It follows that there is a canonical surjective $W^\ast$-homomorphism from $\calr_0(\mathbb R^2)$ onto $R(Q)$. This shows how to apply the above theory even when there are no $C^\ast$-semi-norm (on $Q$).

Similar considerations apply to self-adjoint quantum fields. This leads to the net $\calo\mapsto \calr(\calo)=\calr_0(\cald(\calo))$ of $W^\ast$-algebras where $\cald(\calo)$ is the Schwartz space of $C^\infty$-functions with compact supports in $\calo\subset R^4$  \cite{tre:1967}, \cite{sch:1966}. For developments in this direction see for instance \cite{bor-yng:1992}.


\end{document}